\documentstyle[11pt]{article}
\newtheorem {th}{Theorem}[section]
\newtheorem {lem}[th]{Lemma}

\newtheorem{defn}{Definition}

\def\Cox{\hfill \Box}

\def\sf{\sigma\mbox{-field}}

\def\S{{\cal S}}

\def\E{{\bf{E}}}

\def\P{{\bf{P}}}

\def\Z{{\bf{Z}}}
\def\F{{\cal{F}}}
\def\St{{\tilde{S}}}
\def\Rt{{\tilde{R}}}
\def\ptil{{\tilde{\pi}}}
\def\|{\, | \, }
\def\one{{\bf 1}}

\begin{document}

\begin{titlepage}
\begin{center}
{\large \bf MARKOV CHAINS IN A FIELD OF TRAPS} \\
\end{center}
\vspace{5ex}
\begin{flushright}
Robin Pemantle \footnote{Research supported in part by 
a Sloan Foundation Fellowship, and by a Presidential Faculty
Fellowship}$^,$\footnote{Department of Mathematics, University of
Wisconsin-Madison, Van Vleck Hall, 480 Lincoln Drive, Madison, WI 53706}
  ~\\
Stanislav Volkov \footnote{Laboratory of Large Random Systems, Faculty
of Mathematics and Mechanics, Moscow State University, Russia}
\end{flushright}

\vfill

{\bf ABSTRACT:} \break
We consider a Markov chain on a countable state space, on which is
placed a random field of traps, and ask whether the chain gets trapped
almost surely.  We show that the quenched problem (when the traps are
fixed) is equivalent to the annealed problem (when the traps are updated
each unit of time) and give a criterion for almost sure trapping versus
positive probability of non-trapping.  The hypotheses on the Markov chain
are minimal, and in particular, our results encompass the results of
den Hollander, Menshikov and Volkov (1995).
\vfill


\end{titlepage}

\section{Introduction and statement of results}

\setcounter{equation}{0}

Let $\S$ be a countable space and let $p : \S^2 \rightarrow [0,1]$ be
a set of transition probabilities on $\S$, i.e., $\sum_y p(x,y) = 1$
for all $x \in \S$.  Let $\P_x$ denote any probability measure
such that the random sequence $X_0 , X_1 , X_2 , \ldots$ is a
Markov chain with transition probabilities $\{ p (x,y) \}$,
starting from $X_0 = x$.  Given a function $q : \S \rightarrow [0,1)$,
representing a set of trapping probabilities, we define two new
Markov chains as follows.

\noindent{\bf Quenched problem.}  ``site $x$ is a trap with probability
$q(x)$ forever''.  Let $T \subseteq \S$ be the random set of traps,
i.e., $\P_x ( \{ x_1 , \ldots , x_n \} \subset T) = \prod_{i=1}^n q(x_i)$
for all finite subsets of $\S$, and $T$ is independent of $X_0 , X_1 ,
\ldots$.  We are interested in the quantities
\begin{equation} \label{eq 1}
\pi (x) = \P_x (X_n \in T \mbox{ for some } n \geq 0 ) .
\end{equation}
We say that the quenched field is {\em trapping} or {\em non-trapping},
according to whether $\pi (x) = 1$ for all $x$, or whether $\pi (x) < 1$
for some $x$.

\noindent{\bf Annealed problem.}  ``The state of site $x$ is a trap
with probability $q(x)$, but each unit of time the state of site $x$ is
chosen afresh with this probability''.  This gives us an IID sequence
of trap sets $\{ T_n : n \geq 0 \}$, each distributed as $T$ in the
quenched problem.  We let
\begin{equation} \label{eq 2}
\ptil (x) = \P_x (X_n \in T_n \mbox{ for some } n \geq 0 ) .
\end{equation}
We say that the annealed field is trapping or non-trapping, according
to whether $\ptil (x) = 1$ for all $x$, or whether $\ptil (x) < 1$
for some $x$.

Our first result is the equivalence of these two problems.  Let
$$g(x ,y) = \sum_{n=0}^\infty \P_x (X_n = y)$$
denote the Greens function.

\begin{th} \label{th 1}
Assume there is a constant $K$ such that
\begin{equation} \label{eq green bound}
g(x,x) \leq K \mbox{ for all } x \in \S .
\end{equation}
Then the annealed field is trapping if and only if the
quenched field is trapping.
\end{th}

\noindent{\em Remarks:}

1.  The necessity of an assumption such as~(\ref{eq green bound}) is
illustrated by Example~2 below.

2.  These two problems were considered by den Hollander, Menshikov and
Volkov (1995) in the special case of mean-zero, finite-variance
random walk on $\Z^d$, $d \geq 3$.  They showed (in their Section~5)
that the quenched and annealed problems were equivalent under an additional
technical assumption.  The present theorem thus removes their
technical hypothesis, as well as showing that the particular setting
of $\Z^d$ is irrelevant.

3.  It is elementary (see the next paragraph) that $\pi (x)
\leq \ptil (x)$ for all $x$.  Thus only one direction of Theorem~\ref{th 1}
is interesting.

Let $\tau (x) = \inf \{ k : X_k = x \}$ denote the first hitting time
of the point $x$ and define the two quantities
\begin{eqnarray} \label{eq s1}
R_n & = & \sum_{i=0}^n - \log (1 - q(X_i)) \one_{\tau (X_i) = i} \; ; \\[2ex]
\Rt_n & = & \sum_{i=0}^n - \log (1 - q(X_i)) .  \nonumber
\end{eqnarray}
It is easy to see that the probability of no trapping in the annealed
chain up to time $n$ is given by
$$\P_x (X_i \notin T_i \mbox{ for all } i \leq n \| X_0 , \ldots , X_n)
   = \exp (- \Rt_n) .$$
Similarly, since the conditional probability of $X_n \notin T$ given
$X_0 , X_1 , \ldots$ and given $\{X_i \notin T : i < n \}$, is
equal to $1 - q(X_n) \one_{\tau (X_n) = n}$, the probability of no
trapping in the quenched chain to time $n$ may be computed as
$$\P_x (X_i \notin T \mbox{ for all } i \leq n \| X_0 , \ldots , X_n)
   = \exp (- R_n) .$$
Thus
\begin{eqnarray} \label{eq 5}
\pi (x) & = & \E [1 - \exp (- R_\infty)] \mbox{ and } \\[1ex]
\ptil (x) & = & E [1 - \exp (- \Rt_\infty)] .  \nonumber
\end{eqnarray}
 From equations~(\ref{eq s1}) and~(\ref{eq 5}), it is evident that
$R_n \leq \Rt_n$, and hence that $\ptil (x) \geq \pi (x)$ for all $x$.

We turn now to the determination of a criterion for trapping.
We consider only the annealed case, which of course solves
the quenched case as well under condition~(\ref{eq green bound}).
By~(\ref{eq 5}), the problem reduces to the determination of when
$\Rt_\infty = \infty$ with probability 1.  Define
\begin{eqnarray*}
S_n & = & \sum_{i=0}^n q(X_i) \one_{\tau(X_i) = i} \\[1ex]
\St_n & = & \sum_{i=0}^n q(X_i)
\end{eqnarray*}
and observe that $\Rt_\infty = \infty$ if and only if $\St = \infty$,
and that $R_\infty = \infty$ if and only if $S_\infty = \infty$.
Thus an obvious necessary condition for trapping in the anealed
case is that for all $x_0$,
\begin{equation} \label{eq ss}
\E \St_\infty = \sum_{x \in \S} g(x_0,x) q(x) = \infty .
\end{equation}
Later, we will discuss some conditions under which this is
sufficient as well, but the next example shows that~(\ref{eq ss})
is not always sufficient for trapping.

\begin{defn}
Say that a subset $A \subseteq \S$ is {\em transient for x} if
$\P_x (X_n \in A \mbox{ for some } n) < 1$.  Say that $A$ is
{\em transient} if it is transient for some $x$.  We remark
that when all states communicate in $A^c$, then $A$ must be transient 
for all $x$ or for no $x$; with the present definition, our results 
make sense for more general chains.
\end{defn}

{\bf Example~1}:  Let $A$ be any transient set satisfying
$\sum_{x \in A} g(x_0,x) = \infty$. An example of such a set $A$
for simple random walk in $\Z^3$ is the set
$$\{ x : ||x - (0,0,2^n)|| \leq 2^{n/2} \mbox{ for some } n \geq 1 \} .$$
Let $q(x) = c \in (0,1)$ for $x \in A$ and 0 for $x \notin A$.
Then clearly $\ptil (x_0) < 1$, while $\sum_{x \in \S} g(0,x) q(x)
= \infty$.

Our second result is that this is the only way that~(\ref{eq ss})
can fail to imply trapping.

\begin{th} \label{th 2}
Suppose that the annealed field is non-trapping, i.e., for some $x_0$,
$\ptil (x_0) < 1$.  Then there exists a subset $A \subseteq \S$ such that
\begin{eqnarray*}
(i) && A \mbox{ is transient for } x_0, \hspace{1in} \mbox{ and } \\[1ex]
(ii) && \sum_{x \in \S \setminus A} g(x_0 , x) q(x) < \infty.
\end{eqnarray*}
\end{th}

These two theorems are proved in Section~2.  The final section discusses
conditions under which the condition $\sum_{x \in \S} g(x_0,x) q(x)
= \infty$ is necessary and sufficient for trapping.  These cases
include the case of a simple random walk on $\Z^d$, $d \geq 3$, and
a spherically symmetric function $q(x) = q(||x||)$ discussed in
den Hollander, Menshikov and Volkov (1995).

\section{Proofs of the two theorems}

\setcounter{equation}{0}

\noindent{\sc Proof of Theorem}~\ref{th 1}: By~(\ref{eq 5}),
it suffices to show that the event $\{ S_\infty < \infty = \St_\infty \}$
has probability zero.  Recall that $\tau (x)$ is defined to be the
first hitting time of $x$, and let $\F_{\tau (x)}$ be the $\sf$
generated by the Markov chain up to time $\tau (x)$.  Observe that
for any $x$ and any event $G \in \F_{\tau (x)}$,
\begin{equation} \label{eq 2.1}
\E_{x_0} \one_G \one_{\tau (x) < \infty} ( \sum_{i=0}^\infty
   \one_{X_i = x} )
   = g(x,x) \P_{x_0} (G \cap \{ \tau (x) < \infty \}) \, .
\end{equation}
Given $M \geq 0$, define a time $\tau_M = \inf \{ k : S_k \geq M \}$
and define
$$\St^{(M)} = \sum_{i=0}^\infty q(X_i) \one_{\tau (X_i) \leq \tau_M} .$$

 From the definitions, it is evident that $\St^{(M)} = \St_\infty$
whenever $\tau_M = \infty$.  A second useful fact, also immediate
from the definitions, is that the event $\{ \tau (x) \leq \tau_M \}$
is the same as the intersection of events
$$\{ \tau (x) < \infty \} \cap \{ S_{\tau (x) - 1} < M \} ,$$
from which also we see it is in $\F_{\tau (x)}$.  From these two facts we get
\begin{eqnarray*}
\E_{x_0} \St^{(M)} & = & \sum_{i = 0}^\infty \sum_{x \in \S}
   q(x) \E_{x_0} \one_{X_i = x} \one_{\tau (x) \leq \tau_M} \\[1ex]
& = & \sum_{x \in \S} q(x) \E \sum_{i=0}^\infty \one_{X_i = x}
   \one_{\tau (x) < \infty} \one_{S_{\tau (x) - 1} < M} .
\end{eqnarray*}
Applying~(\ref{eq 2.1}) for each $x$ shows that this is equal to
$$\sum_{x \in \S} q(x) g(x,x) \P_{x_0} (\tau (x) < \infty ,
   S_{\tau (x) - 1} < M) \, ,$$
which is bounded above by
$$\sum_{x \in \S} K q(x) \P_{x_0} (\tau (x) \leq \tau_M)
   = K \E_{x_0} S_{\tau_M} \, .$$
Since $S_{\tau_M} < M+1$, this shows that for any $M$, $\St^{(M)}$
has finite expectation and is hence almost surely finite.
Since $\St^{(M)} = \St_\infty$ whenever $\tau_M = \infty$,
we see that $\St_\infty < \infty$ whenever $\tau_M = \infty$
for some $M$, which happens exactly when $S_\infty < \infty$.  $\Cox$

Before proving Theorem~\ref{th 2}, we give an example to show that
the hypothesis of a bounded Greens function is necessary.

{\bf Example~2}: Let $\S = \{ 2 , 3 , 4  \ldots \}$ with
transition probabilities $p(n,n+1) = 1/n$, $p(n,n) = 1 - 1/n$
and $p(n,k) = 0$ for $k \neq n , n+1$.  The quenched field traps
the chain at $n$ with probability $q(n)$, conditional on its not being
trapped at any $k < n$, hence the quenched field is
trapping if and only if $\sum_n q(n) = \infty$.  The annealed field
traps the chain at $n$ with probability $nq / (nq + 1 - q)$, conditional
on its not being trapped at any $k < n$, hence the annealed field
is trapping if and only if $\sum_n n q(n) = \infty$, which is the
Greens function criterion.

One can replace this example if desired by an example where the Markov
chain is a simple random walk.  Let $G$ be a binary tree, rooted
at a vertex 0, to which is appended at each vertex $v$ in
generation $n > 0$ a single chain of vertices of length $n$.
Setting $q(x) = 1/n$ if $x$ is the end of one of these chains
for some $n$, and $q(x) = 0$ otherwise, makes the annealed field
trapping and the quenched field non-trapping for simple random walk.

We now give the proof of Theorem~\ref{th 2}.  The idea is that
$\ptil (X_n)$ is a martingale if one stops at the value 1 upon being
trapped, and therefore that $\ptil (X_n)$ is a $\P_x$-supermartingale.
The excess of a positive superharmonic function ($\ptil$) is
integrable against the Greens function; this excess is usually
close to $q$, and in fact the set where the excess of $\ptil$ is
not near $q$ must be a transient set.  We begin by stating two
facts whose proofs are immediate.

\begin{lem} \label{lem 1}
For all $x$, under no assumptions on the Markov chain,
$$\ptil (x) = q(x) + (1 - q(x)) \E_x \ptil (X_1) .$$
$\Cox$
\end{lem}

\noindent{\sc Proof of Theorem}~\ref{th 2}:
Let $B_a = \{ x : \ptil (x) \geq a \}$.  Observe that if 
\begin{equation} \label{eq n1}
x \notin B_{1/16} \; \Rightarrow \; P_x (X_1 \in B_{1/4}) < {1 \over 4}.  
\end{equation}
For any $a$, the set $B_a$ is the union of a finite set and a set
transient for $x_0$.  To see this, note that $\ptil (X_n) \rightarrow 0$
almost surely on the event of nontrapping, and hence that
$\P_{x_0} (\ptil (X_n) \rightarrow 0) > 0$; choosing $N$ so that
$\P_{x_0} (\ptil (X_n) < a \forall n \geq N) > 0$, and letting
$F = \{ x_1 , x_2 , \ldots , x_N \}$ such that 
$$\P_{x_0} (\ptil (X_n) < a \forall n \geq N \| X_i = x_i \mbox{ for } 
   i = 1 , \ldots N) > 0,$$ 
gives a set $B_a \setminus F$ that is transient for $x_0$.
We now show that
$\sum_{x \in \S \setminus B_{1/16}} g(x_0 , x) q(x) < \infty$, which
proves the theorem with $A = B_{1/16} \setminus F$. 
By~(\ref{eq n1}), for any $x \notin B_{1/16}$,
$$\E_x \ptil(X_1) \leq {1 \over 4} \P_x (X_1 \notin B_{1/4})
   + \P_x (X_1 \in B_{1/4}) \leq {1 \over 2} \, .$$
Applying Lemma~\ref{lem 1} we see that for any $x$,
\begin{eqnarray*}
\E_x \ptil (X_1) & = & \ptil (x) - q(x) (1 - \E_x \ptil (X_1)) \\[1ex]
& \leq & \ptil (x) - {1 \over 2} q(x) \one_{x \notin B_{1/16}} \, . \nonumber
\end{eqnarray*}

Iterating this, using the Markov property, gives
$$ 0 \leq \E_x \ptil (X_n) \leq \ptil (x) - \sum_{i=0}^{n-1}
   {1 \over 2} \E_x q(X_i) \one_{X_i \notin B_{1/16}}$$
and hence
$$\sum_{i=0}^\infty \E_{x_0} q(x_i) \one_{X_i \notin B_{1/16}} \leq
   2 \ptil (x_0) < \infty .$$
But
$$\sum_{i=0}^\infty \E_{x_0} q(x_i) \one_{X_i \notin B_{1/16}}
   = \sum_{x \in\S \setminus B_{1/16}} g(x_0 , x) q(x) \, ,$$
so this concludes the proof.   $\Cox$.

\section{Chains with well behaved Greens function geometry}

\setcounter{equation}{0}

We have not yet imposed any geometry on $\S$.  In order to
formulate regularity conditions under which~(\ref{eq ss})
is equivalent to trapping, the geometry inherent in the
Greens function must be reasonable, and must be compatible with
the geometry imposed by $q$.  Assume throughout this section
that the Greens function bound~(\ref{eq green bound}) holds.
\begin{defn}
For $x_0 \in \S$, $L \leq K$ and $\alpha \in (0,1)$, define
the Greens function annulus $H_\alpha (L , x_0)$ to be the set
$$\{ x \in \S : L \geq g(x_0 , x) \geq \alpha L \} \, .$$
Say that the Markov chain has {\em reasonable annuli} if
for some $\alpha \in (0,1)$ and for every $L \leq K, \, x_0 \in \S$ and
$A \subseteq \S$ transient with respect to $x_0$, the annulus
$H_\alpha (L , x_0)$ has finite cardinality and
\begin{equation} \label{eq limsup}
\limsup_{L \rightarrow 0} { |H_\alpha (L , x_0) \cap A|
   \over |H_\alpha (L , x_0)|} < 1.
\end{equation}
\end{defn}

Most nice chains satisfy this definition.  For example, consider a
simple random walk in $\Z^d$, $d \geq 3$.  Annuli for this chain
are spherical shells, and any set $A$ taking up more than a fraction
$\beta$ of such a shell is hit with probability at least $\beta + o(1)$.
Thus the probability of hitting $A$ infinitely often is at least the
limsup in~(\ref{eq limsup}), which is less than 1 for a transient
set (actually 0, by tail triviality).  For another example, take the
Markov chain on an infinite rooted binary tree which always walks
away from the root, choosing either of the two children with equal
probability.  This chain is far from irreducible, yet when $\alpha$ is
large enough so that the annuli are nonempty, it clearly
satisfies~(\ref{eq limsup}).

\begin{th} \label{th 3}
Suppose the Markov chain on $\S$ with transitions $p(x,y)$ has
reasonable annuli, and suppose that, for some $C , C' > 1$, the function
$q : \S \rightarrow [0,1)$ satisfies the following regularity condition:
\begin{equation} \label{eq reg}
{1 \over C} g(x_0 , x) \leq g(x_0 , y) \leq C g(x_0 , x) \Rightarrow
   {1 \over C'} q(x) \leq q(y) \leq C' q(x) .
\end{equation}
Then $\ptil (x_0) = 1$ if and only if $\sum_x g(x_0 , x) q(x) = \infty$.
\end{th}

\noindent{\sc Proof:} If $\ptil (x_0) < 1$, then Theorem~\ref{th 2}
gives a set $A$, transient for $x_0$, with
$$\sum_{x \in \S \setminus A} g(x_0 , x) q(x) < \infty \, .$$
By the assumption of reasonable annuli, the transience of $A$ and
condition~(\ref{eq reg}),
$$\limsup_{n \rightarrow \infty} { {\displaystyle \sum_{x \in H_\alpha
   (L \alpha^n , x_0) \cap A} \; \; g(x_0 , x) q(x) } \over
   {\displaystyle \sum_{x \in H_\alpha (L \alpha^n , x_0) \setminus A}
   \; \; g(x_0 , x) q(x)} } < \infty \, ,$$
and hence $\sum_x g(x_0 , x) q(x) < \infty$.   $\Cox$

\noindent{\em Remarks:}

4.  In the case of a simple random walk on $\Z^d$, we find that
condition~(\ref{eq reg}) is equivalent to the existence of a $C < 1$ 
for which
$$\sup_x \sup_{y : C^{-1} |x| \leq |y| \leq C |x|}
   {q(x) \over q(y)} < \infty .$$
This is a more natural formulation than the equivalent condition~(4.1)
in den Hollander, Menshikov and Volkov (1995), namely that
$$\sup_x \sup_{y : |y| \leq \sqrt{d} |x|} {q(x) \over q(y)} < \infty .$$
Hence our Theorem~\ref{th 3} generalizes their integrability test at the
end of their Section~4.1.

5.  Condition~(\ref{eq limsup}) may appear cumbersome but it is more
natural than it looks.  First, it is a condition on the Markov chain
alone, so that Theorem~\ref{th 3} identifies a class of chains such
that the $\sum g(x_0 , x) q(x) = \infty$ is a sharp criterion for
all functions $q$ that are ``spherically symmetric up to constant factor''
as defined by~(\ref{eq reg}).  Secondly, it is nearly sharp, meaning
that if 
$$\lim_{n \rightarrow \infty} { |H_\alpha (L \alpha^n , x_0) \cap A|
   \over |H_\alpha (L \alpha^n , x_0)| } = 1 , $$
then one can always choose a sequence $\{ a_n \}$ such that 
setting $q(x) = a_n$ on $H_\alpha (L \alpha^n , x_0)$ gives
a non-trapping field for which $\sum g(x_0 , x) q(x) = \infty$.

\end{document}